\documentclass{amsart}

 \newtheorem{theorem}{Theorem}[section]
\newtheorem{definition}[theorem]{Definition}

\newtheorem{lemma}[theorem]{Lemma}

\newcommand{\ZZ}{{\mathbb Z}}
\newcommand{\RR}{{\mathbb R}}
\newcommand{\CC}{{\mathbb C}}
\newcommand{\TT}{{\mathbb T}}
\newcommand{\EE}{{\mathbb E}}
\newcommand{\DD}{{\mathbb D}}
\newcommand{\KK}{{\mathbb K}}

\newcommand{\hJ}{\widehat{J}}
\newcommand{\hS}{\widehat{S}}

\newcommand{\nLp}{{L^2_{d\chi}(l^2(\ZZ_+))}}

\title[  Finite difference operators with a finite--band spectrum]{
  Finite difference operators with a finite--band spectrum}

\author{M. Shapiro, V. Vinnikov and P. Yuditskii$^1$}
\thanks{$^1$ This work was supported by the Austrian Science
Found FWF, project number: P16390--N04}

\address{Department of Mathematics,
Michigan State University, East Lansing, MI 48824}
\email{mshapiro@math.msu.edu}
\address{Department of Mathematics,
Ben--Gurion University, P.O. Box 653, Beer--Sheva, 84105, Israel}
\email{vinnikov@cs.bgu.ac.il}
\address{Institute for Analysis,
Johannes Kepler University of Linz, A--4040 Linz, Austria}
 \email{yuditski@math.msu.edu}

\begin{document}

\maketitle

\section{Ergodic finite difference operators and
associated Riemann surfaces}

The standard (three--diagonal) finite--band Jacobi matrices
\cite{DKN, Tes} can be defined as almost periodic or even ergodic
Jacobi matrices with absolutely continuous spectrum that consists
of a finite system of intervals. We wish to find a natural
extension  of this class of finite difference operators onto
multi-diagonal case. First let us recall what is ergodic operator
\cite{Cy, PF}, see also \cite{VYu}.

Let $(\Omega,\frak A, d\chi)$ be a separable
probability space and let
$T:\Omega\to\Omega$ be an invertible
ergodic transformation, i.e., $T$ is
measurable, it preserves $d\chi$, and every
measurable $T$--invariant set has measure
$0$ or $1$.
Let $\{q^{(k)}\}_{k=0}^d$
be functions from $L^\infty_{d\chi}$, more
$q^{(d)}$ is positive--valued and
$q^{(0)}$ is real--valued.

Then with almost every $\omega\in\Omega$ we associate
self--adjoint $2d+1$--diagonal operator $J(\omega)$ as follows:

\begin{equation}\label{general}
(J(\omega)x)_n= \sum_{k=-d}^{d}
\overline{q^{(k)}_{n}(\omega)}x_{n+k}, \quad
          x=\{x_n\}_{n=-\infty}^{\infty} \in l^2(\ZZ),
\end{equation}
where
  $q_n^{(k)}(\omega):=q^{(k)}(T^n\omega)$
and $q^{(-k)}(\omega):= \overline{q^{(k)}(T^{-k}\omega)}$.

Note that the structure of $J(\omega)$ is described by the
following identity
\begin{equation}\label{com}
J(\omega)S=S J(T\omega),
\end{equation}
where $S$ is the shift operator in $l^2(\ZZ)$. The last relation
indicates strongly that one can associate with the family of
matrices $\{J(\omega)\}_{\omega\in\Omega}$ a natural pair of
commuting operators.

Namely, let $L^2_{d\chi}(l^2(\ZZ))$ be the space of
$l^2(\ZZ)$--valued vector functions, $x(\omega)\in l^2(\ZZ)$, with
the norm
\begin{equation*}
||x||^2=\int_\Omega||x(\omega)||^2\,d\chi.
\end{equation*}
Define
\begin{equation*}
(\hJ x)(\omega)=J(\omega)x(\omega), \quad (\hS x)(\omega)=S
x(T\omega), \quad x\in L^2_{d\chi}(l^2(\ZZ)).
\end{equation*}
Then \eqref{com} implies
\begin{equation*}
(\hJ\hS x)(\omega)= J(\omega)S x(T\omega)= S J(T\omega)
x(T\omega)= (\hS\hJ x)(\omega).
\end{equation*}

Further, $\hS$ is a unitary operator and
$\hJ$ is self-adjoint. The space
$L^2_{d\chi}(l^2(\ZZ_+))$
is an invariant subspace for $\hS$. It is not
invariant with respect to $\hJ$ but it is invariant with respect to
$\hJ\hS^d$.
Put
\begin{equation*}
\hS_+=\hS| L^2_{d\chi}(l^2(\ZZ_+)),
\quad
(\hJ\hS^d)_+=\hJ\hS^d| L^2_{d\chi}(l^2(\ZZ_+)).
\end{equation*}


\begin{definition}[local functional model]
We say that a pair of commuting operators $A_1:H\to H$ and
$A_2:H\to H$ has a (local) functional model if there is a unitary
embedding $i:H\to H_O$ in a space $H_O$ of holomorphic in some
domain $O$ functions $F(\zeta)$, $\zeta\in O$, with a reproducing
kernel ($F\mapsto F(\zeta_0)$, $\zeta_0\in O$, is a bounded
functional in $H_O$) such that operators $i_\star A_1$ and
$i_\star A_2$ become a pair of operators of multiplication by
holomorphic functions, say
\begin{equation*}
A_1 x\mapsto a_1(\zeta) F(\zeta), \quad A_2 x\mapsto a_2(\zeta)
F(\zeta).
\end{equation*}
\end{definition}

Existence of a local functional model implies a number of quite
strong consequences. In what follows $b(\zeta)$ and $\lambda
(\zeta)$ denote the functions (symbols) related to operators
$\hS_+$ and $(\hJ\hS^d)_+$. Let $k_\zeta$ be the reproducing
kernel in $H_O$ and let $\hat k_\zeta$ be its preimage $i^{-1}
k_\zeta$ in $L^2_{d\chi}(l^2(\ZZ_+))$. Then
\begin{equation}\label{3}
\langle \hS_+^* \hat k_\zeta, x\rangle
=\langle \hat k_\zeta, \hS_+ x\rangle=
\langle k_\zeta, b F\rangle.
\end{equation}
By the reproducing property
\begin{equation*}
\langle k_\zeta, b F\rangle=
\overline{b(\zeta) F(\zeta)}=
\langle \overline{b(\zeta)}k_\zeta,  F\rangle.
\end{equation*}

Hence,
\begin{equation*}
\langle \hS_+^* \hat k_\zeta, x\rangle=
\langle \overline{b(\zeta)}\hat k_\zeta,
x\rangle.
\end{equation*}
That is $\hat k_\zeta$ is an eigenvector of $\hS^*_+$ with the
eigenvalue $\overline {b(\zeta)}$. In the same way, $\hat k_\zeta$
is an eigenvector of $(\hJ\hS)^*_+$ with the eigenvalue
$\overline{\lambda(\zeta)}$.

Thus, if a functional model exists then
the spectral problem
\begin{equation}\label{4}
\left\{\begin{matrix}
  \hS_+^* \hat k_\zeta &=
\overline{b(\zeta)}\hat k_\zeta\\
(\hJ\hS^d)_+^* \hat k_\zeta &=
\overline{\lambda(\zeta)}\hat k_\zeta
\end{matrix},\right.
\end{equation}
has a solution $\hat k_\zeta $
antiholomorphic in $\zeta$.
 Moreover, linear combinations of all $\hat
k_\zeta$ are dense in $L^2_{d\chi}(l^2(\ZZ_+))$. Vice versa, if
\eqref{4} has a solution of such kind  then we define
\begin{equation*}
F(\zeta):=\langle x, \hat k_\zeta\rangle,
\quad ||F||^2:=||x||^2.
\end{equation*}
This provides a local functional
model for the pair  $\hS_+$,
  $(\hJ\hS)_+$.

The following lemma is evident.
\begin{lemma} Let $U:L^2_{d\chi}\to
L^2_{d\chi}$ be the unitary operator associated with the ergodic
transformation $T$: $(Uc)(\omega)=c(T\omega)$,
$c\in L^2_{d\chi}$.
We denote by the
same letter $q$ both a function $q\in L^\infty_{d\chi}$ and the
 multiplication
operator by $q$ (e.g., $(q c)(\omega):= q (\omega) c(\omega)$).

Problem \eqref{4} is equivalent to the following spectral problem
\begin{equation}\label{5}
\left\{
\sum_{k=-d}^d U^k
\overline{q^{(k)}b^k(\zeta)}\right\}
c_\zeta=\overline{z(\zeta)}c_\zeta,
\end{equation}
where $z(\zeta):=\lambda(\zeta)/b^d(\zeta)$ and $c_\zeta$ is an
anti--holomorphic $L^2_{d\chi}$-- valued vector function. Moreover
$\{c_\zeta\}$ is complete in $L^2_{d\chi}$ if and only if $\{\hat
k_\zeta\}$ is complete in $L^2_{d\chi}(l^2(\ZZ_+))$.
\end{lemma}

We may hope to glue a {\it global} functional model on a Riemann
surface $X_0=\DD/\Gamma_0$ formed by functions $(z,b)$. Of course,
this model does not necessarily exist (even existence of a local
model requires some additional assumptions on the ergodic map and
the coefficients functions).

The surface $X_0$ is in generic case of infinite genus. However we can
reduce it because $X_0$ possesses a family of automorphisms.
Let $e_{\{\gamma\}}$ be an eigenvector of $U$ with an eigenvalue
$\bar\mu_{\{\gamma\}}$.
The systems of eigenfunctions and eigenvalues form both Abelian groups with
respect to multiplication.
Using \eqref{5} we get immediately that $\{\gamma\}:(z,b)\mapsto
(z,\mu_{\{\gamma\}}b)$ is an automorphism of $X_0$. Taking a
quotient of $X_0$ with respect to these automorphisms we obtain a
much smaller surface $X=\DD/\Gamma$, $\Gamma=\{\gamma\Gamma_0\}$.
  Note that $z$ is still a function on $X$ but
$b$ becomes a character automorphic function. Finally, using $z$
we may glue boundary of $X$ and to get in this way a compact
 Riemann surface $X_c$, such that  $X=X_c\setminus E$. The simplest
assumption is that the boundary $E$ is a finite system of cuts on
$X_c$.

Thus in a certain case we may expect that the  triple
$\{X_c,z,E\}$ characterizes spectrum  of a finite difference
operator. If so we call it a finite band operator.  We point out
that in fact every triple of this kind gives rise to a family of
ergodic finite difference operators. More precise and detailed
characteristic of such operators is given in the next section.

Now we would like to describe all triples of given type up to a
natural equivalence relation.

\begin{definition}
We say that two triples $({X_c}_1,z_1,E_1)$
and $({X_c}_2,z_2,E_2)$ are equivalent
if there exists a holomorphic homeomorphism
$h:{X_c}_1\to {X_c}_2$ such that
$z_1=h^*(z_2)$ and $E_2=h(E_1)$.
\end{definition}

Note, that for any triple $(X_c,z,E)$ the
holomorphic function $z:X_c\to\CC P^1$ is a
ramified coverings of $\CC P^1$. It is quite
convenient to describe equivalence classes of
ramified coverings in terms of branching
 divisor. Namely, a point $p\in X_c$ such that
$\frac{dz}{d \zeta}|_p=0$ where $\zeta$ is a
local holomorphic coordinate in a
neighborhood of $p$ is called
\emph{ramification point (or, critical
point)}. Its image $z(p)$ is called \emph{a
branching point (or, critical value)}. The set
of all branching points of function $z$ form a
branching divisor $B(z)$ of $z$.

Clearly, branching divisors of
equivalent functions are the same. Moreover,
the compact holomorphic curve $X_c$ is also
uniquely determined by the branching divisor
and some additional ramification data (of
combinatorial type). Namely, assume that $z$
has degree $d$ and  let $z_i\in\CC
P^1,\,\,i=1,\dots,N$ be branching points of
function $z$,
$w\in \CC P^1$ be a non branching point.
Fix a system of non-intersecting paths
$\gamma=\{\gamma_i\}_{i=1,\dots,N}$. The $i$th
path $\gamma_i$ connects $w$ and $z_i$. We
want to construct a system of loops
$l_i\subset \CC P^1$. To construct $l_i$ we
start from $w$ and follow first $\gamma_i$
almost to $z_i$, then encircle
$z_i$ counterclockwise along a small circle and finally go back to $w$
along $-\gamma_i$.
Using $l_i$ we associate with each branching point an element of
permutation group
$\sigma_i\in \Sigma_d$.

The point $w$ has exactly $d$ preimages.
Let us label them by integers $\{1,\dots,d\}$.
Let us follow the loop $l_i$ and lift this loop to $X_c$ starting
from each of the preimages
of $w$.
The monodromy along path $l_i$ gives us a permutation
$\sigma_i\in\Sigma_d$ of preimages.
 Therefore, the function $z$ determines $N$ branching points and $N$
permutations.
These permutations are not uniquely defined,
they depend on the labeling of preimages of
$w$. Therefore, they are determined up to a
conjugacy by the elements of $\Sigma_d$.

Given a set of branching
points $B(z)=\{z_i\}_{i=1,\dots,N}\subset\CC
P^1$ and a system of permutations
$\sigma(z)=(\sigma_1,\dots,\sigma_N)\in 
\Sigma_d\times\dots\times \Sigma_d/
\Sigma_d$, where the last quotient is taken
 with respect to diagonal conjugation,
we can restore by Riemann theorem
the surface $X_c$ and function $z$.

Hence, the triple $(X_c,z,E)$ is equivalent
to the quadruple $(B(z),\gamma,\sigma(z),E)$.

\section{The Global Functional Model}

Let $\pi(\zeta):\DD\to X$ be a
uniformization of the surface $X=X_c\setminus
E$.  Thus there exists
a discrete subgroup $\Gamma$ of the group
$SU(1,1)$
consisting of elements of the form
\begin{equation*}
\gamma=\bmatrix \gamma_{11}&\gamma_{12}\\
\gamma_{21}
&\gamma_{22}
\endbmatrix,\
\gamma_{11}=\overline{\gamma_{22}},\
\gamma_{12}=
\overline{\gamma_{21}},
\ \det\gamma=1,
\end{equation*}
such that  $\pi(\zeta)$
is automorphic with respect to $\Gamma$, i.e.,
$\pi(\gamma(\zeta))=\pi(\zeta),\ \forall \gamma
\in\Gamma$,
 and any two preimages of $P_0\in X$ are
$\Gamma$--equivalent.
  We normalize $Z(\zeta):=(z\circ\pi)(\zeta)$ by
the conditions
$Z(0)=\infty$,
$(\zeta^d Z)(0)>0$.

Note that $\Gamma$
acts dissipatively on
$\TT$ with respect to the Lebesgue measure
$dm$, that is there exists a  measurable
(fundamental) set
$\EE$, which does not
contain any two
$\Gamma$--equivalent points,
and the union
$\cup_{\gamma\in\Gamma}\gamma(\EE)$
is a set of
full measure.
In fact $\EE$ can be chosen as a finite union
of intervals. For the space of square
summable functions  on $\EE$ (with respect
to $dm$), we use
the notation $L^2_{dm|\EE}$.

A character of $\Gamma$ is a complex--valued
function
$\alpha:\Gamma\to \TT$, satisfying
\begin{equation*}
\alpha(\gamma_1\gamma_2)=\alpha(\gamma_1)
\alpha(\gamma_2),\quad\gamma_1,\gamma_2
\in\Gamma.
\end{equation*}
The characters form an Abelian compact
group denoted
by $\Gamma^*$.

Let $f$ be an analytic function in
$\DD$, $\gamma\in\Gamma$.
Then we put
\begin{equation*}
f\vert[\gamma]_k=\frac{f(\gamma(\zeta))}
{(\gamma_{21}\zeta+\gamma_{22})^k}
\quad k=1, 2.
\end{equation*}
Notice that $f\vert[\gamma]_2=
f$ for all $\gamma\in \Gamma$, means
that the form
$f(\zeta)d\zeta$ is invariant with respect to
the substitutions
$\zeta\to\gamma(\zeta)$ ($f(\zeta)d\zeta$ is
an Abelian
integral on $\DD/\Gamma$). Analogously,
$f\vert[\gamma]=
\alpha(\gamma)f$ for all $\gamma\in \Gamma$,
$\alpha\in\Gamma^*$,
means
that the form
$|f(\zeta)|^2\,|d\zeta|$ is invariant with
respect to
these substitutions.

We recall, that a function $f(\zeta)$ is
of Smirnov
class, if it can be represented as a ratio
of two
functions from $H^\infty$ with an outer
denominator. The following spaces related
to the Riemann surface $\DD/\Gamma$ are
counterparts of the standard Hardy spaces $H^2$
($H^1$) on the unit disk.

\begin{definition} \label{def}
The space $A^{2}_1(\Gamma,\alpha)$
($A^{1}_2(\Gamma,\alpha)$) is formed by functions $f$,
which are analytic on $\Bbb D$ and satisfy the
following three conditions
\begin{equation*}
\begin{split}
1)& f \ \text{is of Smirnov class}\\
2)& f\vert[\gamma]=\alpha(\gamma) f\ \ \
(f\vert[\gamma]_2=\alpha(\gamma) f)
\quad
\forall\gamma\in \Gamma\\
3)&
\int_{\Bbb E}\vert f\vert^2\,dm<\infty\  \ \
(\int_{\Bbb E}\vert f\vert\,dm<\infty).
\end{split}
\end{equation*}
\end{definition}

$A^2_1(\Gamma,\alpha)$
is a Hilbert space with the reproducing kernel
$k^\alpha(\zeta,\zeta_0)$, moreover
\begin{equation}\label{widom}
0<\inf_{\alpha\in\Gamma^*}
k^\alpha(\zeta_0,\zeta_0)\le
\sup_{\alpha\in\Gamma^*}
k^\alpha(\zeta_0,\zeta_0)
<\infty.
\end{equation}
Put
\begin{equation*}
k^\alpha(\zeta)=k^\alpha(\zeta,0)\quad\text{and} \quad
K^\alpha(\zeta)=\overline{K^\alpha_\zeta(0)}=\frac{k^\alpha(\zeta)}
{\sqrt{k^\alpha(0)}}.
\end{equation*}

We need one more special function.
  The Blaschke product
\begin{equation*}
b(\zeta)=\zeta
\prod_{\gamma\in\Gamma, \gamma\not= 1_2}
\frac{\gamma(0)-\zeta}{1-\overline{\gamma(0)}
\zeta}
  \frac{\vert\gamma(0)\vert}{\gamma(0)}
\end{equation*}
is called the {\it Green's function}
of $\Gamma$ with
respect to the origin.
It is a character--automorphic function,
i.e., there
exists $\mu\in\Gamma^*$ such that
$b(\gamma(\zeta))=\mu(\gamma)
b(\zeta)$.
  Note, if $G(P)=G(P,\infty)$ denotes
the Green's function of the surface
$X$, then
\begin{equation*}
G(\pi(\zeta))=-\log\vert b(\zeta)\vert.
\end{equation*}

Let $\Gamma_0:=\ker\mu$, that is
$\Gamma_0=\{\gamma\in\Gamma:
\mu(\gamma)=1\}$.  Evidently, $b(\zeta)$ and
$(Zb^d)(\zeta)$ are holomorphic functions on
the surface
$X_0=\DD/\Gamma_0$.

Now, assume that $\alpha_0\in\Gamma_0$ can be
extended to a character on $\Gamma$, i.e.,
\begin{equation*}
\Omega_{\alpha_0}=
\{\alpha\in\Gamma^*: \alpha|\Gamma_0=
\alpha_0\}\not=\emptyset.
\end{equation*}
Note that the set of characters
\begin{equation*}
\Omega_{\iota}=
\{\alpha\in\Gamma^*: \alpha|\Gamma_0=
\iota\}
\end{equation*}
where $\iota(\gamma)=1$ for all
$\gamma\in\Gamma_0$ is isomorphic to
the set $(\Gamma/\Gamma_0)^*$.

Let us fix an element
$\hat\alpha_0\in\Omega_{\alpha_0}$. Since
\begin{equation*}
  \{\alpha\in\Gamma^*: \alpha|\Gamma_0=
\alpha_0\} =
\{\hat\alpha_0\beta:\beta\in\Gamma^*:
\beta|\Gamma_0=
\iota\}
\end{equation*}
we can define a measure
$d\chi_{\alpha_0}(\alpha)$ on
$\Omega_{\alpha_0}$ by the relation
\begin{equation*}
d\chi_{\alpha_0}(\alpha)=
d\chi_{\alpha_0}(\hat\alpha_0\beta)
=d\chi_{\iota}(\beta),
\end{equation*}
where $d\chi_{\iota}(\beta)$ is the Haar
measure on $(\Gamma/\Gamma_0)^*$
(the measure $d\chi_{\alpha_0}(\alpha)$
does not depend on a choice of the
element $\hat\alpha_0$).

Obviously,
  $T\alpha:=\mu^{-1}\alpha$ is an invertible
ergodic measure--preserving transformation
on $\Omega=\Omega_{\alpha_0}$
with respect to the measure
$d\chi=d\chi_{\alpha_0}$.

The following Theorem is a slightly modified version of 
Theorem 2.2 from    \cite{VYu}.

\begin{theorem}\label{model}
Given $\alpha\in\Gamma^*$, the system of functions $\{b^n
K^{\alpha\mu^{-n}}\}_{n\in\ZZ}$ forms an orthonormal basis in
$L^2_{dm|\EE}$. With respect to this basis, the multiplication
operator by $Z$ is  a $2d+1$--diagonal ergodic finite difference
operator with $\Omega=\Omega_{\alpha_0}$,
  $d\chi=d\chi_{\alpha_0}$,
  $T\alpha:=\mu^{-1}\alpha$ and
$\alpha_0=\alpha|\Gamma_0$.

Moreover, the operators $\hS_+$ and
$(\hJ\hS^d)_+$ are unitary equivalent to
multiplication by $b$ and $(b^d Z)$
in $A^2_1(\Gamma_0,\alpha_0)$ respectively.
This unitary map is given by the formula
\begin{equation*}
\sum_{\{\gamma\}\in\Gamma/\Gamma_0}
f|[\gamma]\alpha^{-1}(\gamma)=
\sum_{n\in\ZZ_+}x_n(\alpha)b^n
K^{\alpha\mu^{-n}},
\end{equation*}
where
$f\in A_1^2(\Gamma_0, \alpha_0)$
and the vector function
$x(\alpha)=
\{x_n(\alpha)\}$ belongs to $\nLp$.
\end{theorem}

\section{Uniqueness Theorem}
\begin{theorem}
Assume that a finite difference ergodic operator has a finite band
functional model that is there exist a triple $\{X_c, \tilde z, E
 \}$, a character $\alpha_0\in\Gamma^*$ and a map $F$ from $\Omega$
to $\tilde\Omega:=\Omega_{\alpha_0}$ such that
$FT\omega=\mu^{-1}F\omega$,  $\chi (F^{-1}(A))=\tilde\chi(A)$,
$A\subset\tilde \Omega$, with $d\tilde\chi:=d\chi_{\alpha_0}$,
here $\mu$ is the character of the Green's function $\tilde b$ on
$X_c\setminus E$. Moreover $q^{(k)}(\omega)=\tilde
q^{(k)}(F\omega)$, where the coefficients $\tilde q^{(k)}(\alpha)$
are generated by the multiplication operator $\tilde z$ with
respect to the orthonormal basis $\{\tilde b^n
K^{\alpha\mu^{-n}}\}_{n\in\ZZ}$.

If the functions $\tilde z$ and $\{d\log\tilde b/d\tilde z\}$
separate points on $X_c\setminus E$ then any local functional
model is generated by one of the branches of the function $\tilde
b$.
\end{theorem}
\begin{proof}

Put
\begin{equation}
A(z;\omega)=
\begin{bmatrix}
0&1& &\dots&0\\
\vdots& &\ddots& &\vdots\\
\vdots& & &\ddots &\vdots\\
0&0& &\dots&1\\
-\overline{\left(\frac{q^{(-d)}(T^{-d}\omega)}
{q^{(d)}(T^{d}\omega)}\right)}&\dots
&\overline{\left(\frac{z(\zeta)-q^{(0)}(\omega)}
{q^{(d)}(T^{d}\omega)}\right)}&\dots
&-\overline{\left(\frac{q^{(d-1)}(T^{d-1}\omega)}
{q^{(d)}(T^{d}\omega)}\right)}
\end{bmatrix}.
\end{equation}

According to \eqref{5}
\begin{equation}\label{8}
A(z;\omega) f_\zeta(\omega) =\overline{b(\zeta)} f_\zeta(T\omega)
\end{equation}
with
\begin{equation*}
f_\zeta(\omega)=
\begin{bmatrix}
c_\zeta(T^{-d}\omega)\overline{b^{-d}(\zeta)}\\
\vdots\\
 c_\zeta(T^{d-1}\omega)\overline{b^{d-1}(\zeta)}
\end{bmatrix}.
\end{equation*}

Considering (if necessary) a subdomain $\tilde O\subset O$ let us
introduce one--to--one map $\zeta\mapsto z$. Then, let us put
$\tilde z= z$ assuming that $\{\tilde z, d\log\tilde b/d\tilde
z\}$ is a point on $X_c\setminus E$. Actually there are exactly
$d$ different preimages $\{\tilde b_l\}_1^d$ with the given
 $\tilde z$
and different $d\log\tilde b_l/d\tilde z$. Finally, since a pair
$(\tilde z,\tilde b_l)$ determines a point on $X_c\setminus E$ we
can put $K(\alpha;\tilde z,\tilde b_l):=K_{\tilde\zeta}^\alpha(0)$
choosing one of preimages $\tilde\zeta$ on the universal covering.

 Using this notation  and identities $\tilde z= z$, 
 $q^{(k)}(\omega)=\tilde
q^{(k)}(F\omega)$, we write
\begin{equation}\label{9}
A(z;\omega) \KK(z;F\omega) = \KK(z;\mu^{-1}F\omega)
\begin{bmatrix}
\overline{\tilde b_1}& & \\
 &\ddots& \\
 & &\overline{\tilde b_d^{-1}}
\end{bmatrix},
\end{equation}
where the matrix $\KK$ is constructed from the reproducing kernels
$$
\KK(z,\alpha)=
\begin{bmatrix}
K(\mu^{d}\alpha;z,\tilde b_1)\overline{\tilde b_1^{-d}}& \dots&
K(\mu^{d}\alpha;z,\tilde b_d^{-1})\overline{\tilde b_d^{d}}
\\
\vdots& &\vdots \\
K(\mu^{-d+1}\alpha;z,\tilde b_1)\overline{\tilde b_1^{-d}}& \dots&
K(\mu^{-d+1}\alpha;z,\tilde b_d^{-1})\overline{\tilde b_d^{d}}
\end{bmatrix}.
$$

Combining \eqref{8} and \eqref{9} we get
$$
\begin{bmatrix} \overline{\tilde b_1(\zeta)}h_1(\omega)\\ \vdots\\
\overline{\tilde b_d^{-1}(\zeta)}h_{2d}(\omega)
\end{bmatrix}=\overline{b(\zeta)}
\begin{bmatrix} h_1(T\omega)\\ \vdots\\
h_{2d}(T\omega)
\end{bmatrix}
$$
with
$$
\begin{bmatrix} h_1(\omega)\\ \vdots\\
h_{2d}(\omega)
\end{bmatrix}:=\KK(z,\alpha)^{-1}f_\zeta(\omega).
$$

First of all, $h_{d+l}(\omega)=0$ for $l\ge 1$ because $|\tilde
b_l(\zeta) b(\zeta)|<1$. Since the spectrum of $U$ is discreet and
$\tilde b_l(\zeta)$ as well as $b(\zeta)$ are holomorphic the
ratio $\tilde b_l(\zeta)/b(\zeta)$ should be a constant if only
$h_l(\omega)\not= 0$. Making use the assumption $(\log\tilde
b_l(\zeta))'\not=(\log\tilde b_m(\zeta))'$, $l\not= m$, we obtain
that only one of entries $h_{l_0}(\omega)$ is different from zero.

Therefore, $
c_\zeta(\omega)=e_{\{\gamma_0\}}(\omega)K(F\omega;z,\tilde
b_{l_0}) $ with
$Ue_{\{\gamma_0\}}=\bar\mu_{\{\gamma_0\}}e_{\{\gamma_0\}}$; and
hence $Af:=e_{\{\gamma_0\}}(\omega)f(F\omega)$ is not only
isometric but a unitary map from $L^2_{d\tilde\chi}$ to
$L^2_{d\chi}$ such that $A\tilde U= UA$. This means that
$\bar\mu_{\{\gamma_0\}}$ is an eigenvalue of $\tilde U$ and the
local model is given by one of the branches of the function
$\tilde b$, $c_\zeta(\omega)=K(F\omega;z,\mu_\gamma\tilde
b_{l_0})$ with a certain $\gamma$.
\end{proof}

The following example shows that in the case when the functions
$\tilde z$ and $\{d\log\tilde b/d\tilde z\}$ do not separate
points on $X_c\setminus E$ one can give different {\it global}
functional realizations for the same ergodic
operator.

\begin{proof}[Example]
Let $J=S^d+S^{-d}$. There exist a "trivial" functional model with
$X_c\setminus E\sim \DD$. In this case $J$ is the multiplication
operator by $z=\zeta^d+\zeta^{-d}$ with respect to the standard
basis $\{\zeta^l\}$ in $L^2_\TT$. Note that $b=\zeta$, thus
$$
w:=\frac{d\log b}{dz}=\frac 1{\zeta^d-\zeta^{-d}}\frac 1 d,
$$
that is $z^2+(w d)^{-2}=4$, $|(wd)^{-1}+z|<2$.

On the other hand let us fix any polynomial $T(u)$, $\deg T=d$,
with real critical values on $\RR\setminus[-2,2]$ and define
$X_c\setminus
E=T^{-1}(\overline\CC\setminus[-2,2])\sim\overline\CC\setminus
T^{-1}[-2,2]$. As it well known the last set is the resolvent set
for a $d$--periodic Jacobi matrix \cite{Tes}, say $J_0$. Moreover
$T(J_0)=J$, and $-\log|b|$ is just the Green's function of this
domain in the complex plain. So, using the standard functional
model for $J_0$ with the symbols $u$ and $b$ we get a functional
model for $J$ with $z=T(u)$ and the same $b$. Note that as before
$z^2+(w d)^{-2}=4$, $|(wd)^{-1}+z|<2$ with $w:=\frac{d\log
b}{dz}$.

\end{proof}

\bibliographystyle{amsplain}

\end{document}